\documentclass{article}      
\usepackage{amssymb}
\usepackage{amsmath}
\usepackage{amscd}
\usepackage{amsthm}  
\usepackage{epsfig}  
\usepackage{subfigure}
\usepackage{graphicx}
\usepackage{amsfonts}
\usepackage{mathrsfs}
\usepackage{units}              % for "nicefrac"
\usepackage{hyperref}
 
\usepackage[all]{xy}
\entrymodifiers={+!!<0pt,\fontdimen22\textfont2>}
\def\objectstyle{\displaystyle}

\usepackage{color}
\definecolor{fixcolor}{rgb}{0.65,0.0,0.00}

\newcommand{\inlinexy}[1]{{\def\objectstyle{\scriptstyle}\def\labelstyle{\scriptstyle}{\xymatrix@-1.2pc{#1}}}}

\usepackage{wrapfig}
\newtheorem{theorem}{Theorem}
\newtheorem{corollary}[theorem]{Corollary}
\newtheorem{lemma}[theorem]{Lemma}  
\newtheorem{proposition}[theorem]{Proposition}

\newtheorem{definition}{Definition}
\newcommand{\bc}{\mathbb{C}}
\newcommand{\bp}{\mathbb{ P}}
\newcommand{\bz}{\mathbb{Z}}
\newcommand{\br}{\mathbb{R}}

\newcommand{\p}{\partial}

\newcommand{\ch}{\mathcal{H}}

\newcommand{\med}{\medskip}
\newcommand{\la}{\longrightarrow}
\newcommand{\bfl}{\begin{flushleft}}
\newcommand{\efl}{\end{flushleft}}

\newcommand{\vn}{\vec{n}}
\newcommand{\vm}{\vec{m}}
\newcommand{\vo}{(0)}

\newcommand{\comp}{\circ}       % Function composition
\newcommand{\R}{{\mathbb R}}        % The reals
\newcommand{\set}[1]{\left\{#1\right\}}
\newcommand{\bdry}{\partial}
\newcommand{\tensor}{\otimes}

\newcommand{\cross}{\times}

\newcommand{\isom}{\cong}

\newcommand{\htpyeq}{\simeq}

\newcommand{\bigdisjcup}{\coprod}
\renewcommand{\emptyset}{\varnothing}
\newcommand{\ob}{\operatorname{Ob}}
\newcommand{\mor}{\operatorname{Mor}}
\newcommand{\nmor}[1]{\operatorname{{\scriptstyle#1}-Mor}}
\newcommand{\diff}{\operatorname{Diff}}
\newcommand{\hocolim}{\underrightarrow{\operatorname{hocolim}}}

\newcommand{\diffocp}{\operatorname{Diff^+_{\textrm{oc}}}}

\newcommand{\defem}[1]{\emph{#1}} % definiton emphasis

\newcommand{\bdryI}{\bdry_I}
\newcommand{\bdryF}{\bdry_{\mathrm{f}}}

\renewcommand{\phi}{\varphi}

\newcommand{\s}{\mathscr{S}}
\newcommand{\scrb}{\mathscr{B}}

\newcommand{\soc}{\mathscr{S}^\textrm{oc}}
\newcommand{\socD}{\mathscr{S}^\textrm{oc}_D}
\newcommand{\socG}{\mathscr{S}^\textrm{oc}_\Gamma}
\newcommand{\socDb}{\mathscr{S}^{\textrm{oc},b}_D}
\newcommand{\socGb}{\mathscr{S}^{\textrm{oc},b}_\Gamma}
\newcommand{\B}{{\mathscr B}}  % topologification functor
\newcommand{\catC}{{\mathscr C}}        % Category C

\newcommand{\tSigma}{\widetilde{\Sigma}}
\begin{document}
\title{The topology of the category of open and closed strings}
 
\author{Nils A. Baas 
  \quad Ralph L. Cohen 
  \thanks{The second author was partially supported by a grant from the NSF} 
  \quad Antonio Ram{\'\i}rez } \date{\today} \maketitle
 \begin{centerline} {\sl Dedicated to Sam Gitler on the occasion of his 70th birthday.}\end{centerline}

%\tableofcontents

\begin{abstract}
  In this paper we study the topology of the cobordism category $\soc$
  of open and closed strings.  This is a 2-category in which the
  objects are compact one-manifolds whose boundary components are
  labeled by an indexing set (the set of ``$D$-branes''), the
  $1$-morphisms are cobordisms of manifolds with boundary, and the
  $2$-morphisms are diffeomorphisms of the surface cobordisms.  Our
  methods and techniques are direct generalizations of those used by
  U.~Tillmann in her study of the category of closed strings. We input
  the striking theorem of Madsen and Weiss regarding the topology of
  the stable mapping class group to identify the homotopy type of the
  geometric realization $|\soc|$ as an infinite loop space.
\end{abstract}

\section*{Introduction}  
As described by Segal \cite{segal}, in a conformal field theory one
associates to a closed $1$-dimensional manifold $S$ a Hilbert space
$\ch_S$, and to a $2$-dimensional cobordism $\Sigma$ between closed
$1$-manifolds $S_1$ and $S_2$, an operator $\phi_\Sigma : \ch_{S_1}
\to \ch_{S_2}$.  This operator depends on the conformal class of a
metric on the surface $\Sigma$.  Such a field theory can be viewed as
a representation of the category $\s$ whose objects are closed
$1$-manifolds, and where the space of morphisms from $S_1$ to $S_2$ is
the moduli space of conformal classes of metrics (i.e Riemann
surfaces) on cobordisms from $S_1$ to $S_2$.
 
In \cite{tillmann} U. Tillmann described $\s$ as a $2$-category, whose
objects are nonnegative integers, corresponding to the number of
components in the closed $1$-manifold, the $1$-morphisms between $m$
and $n$ are cobordisms between a disjoint union of $m$-circles and a
disjoint union of $n$-circles, and whose $2$-morphisms are
diffeomorphisms of the cobordisms that are fixed on the boundary.  She
constructed $\s$ as a symmetric monoidal $2$-category, and with a
clever argument using the Harer stability theorem, she proved that
there is a homotopy equivalence
$$
\bz \times B\Gamma_\infty^+ \simeq \Omega |\s|.
$$
Here $|\s|$ is the bisimplicial geometric realization (classifying
space) of the category $\s,$ and $\Omega |\s|$ is its based loop
space.  $B\Gamma_\infty^+$ is Quillen's plus construction on the
classifying space of the stable mapping class group,
$$
B\Gamma_\infty = \varinjlim_g B\Gamma_{g,1}
$$
where $\Gamma_{g,n}$ is the mapping class group of surfaces of
genus $g$ with $n$ boundary components.  Since $\s$ is a symmetric
monoidal category this result established $B\Gamma_\infty^+$ as an
infinite loop space.
 
This result indicated the important role of infinite loop space theory
in studying the stable mapping class group, and it eventually led to
Madsen and Weiss's proof of a generalization of a conjecture of Mumford
\cite{madweiss}, saying that there is an equivalence of infinite loop
spaces
$$
B\Gamma_\infty^+  \simeq \Omega^\infty \bc \bp^\infty_{-1}.
$$
Here $\Omega^\infty \bc \bp^\infty_{-1}$ is the zero space of the
stunted projective spectrum $\bc \bp^\infty_{-1},$ which is defined to
be the Thom spectrum of minus the canonical line bundle $L \to \bc
\bp^\infty$, $ \bc \bp^\infty_{-1} = (\bc \bp^\infty)^{-L}.$
 
\med The purpose of this paper is to generalize Tillmann's results to
the category $\soc$ of ``open and closed strings.''  In this category
the objects are compact $1$-manifolds which may have boundary.  In
other words, the objects are disjoint unions of circles and closed
intervals.  The need to study such a category arises from open-closed
string theory, in the presence of ``$D$-branes'' (see for example,
\cite{segal2}). If the background space of the string theory is a
manifold $M$, the $D$-branes are submanifolds $B_i \subset M$ that
define the boundary conditions of the strings.  In other words, one is
considering spaces of closed and ``open'' strings moving in $M$, where
the open strings are maps from the interval $\gamma : [0,1] \to M$,
with boundary conditions $\gamma (0) \in B_0$, and $\gamma (1) \in
B_1$, where $B_0$ and $B_0$ are $D$-branes.  In practice, $D$-branes
are often equipped with bundles and connections, but we will not be
concerned with this additional structure in this paper.
 
In order to encode the additional information the presence of
$D$-branes imposes, we consider an abstract finite set $\scrb$ which
is meant to index the set of $D$-Branes in the background manifold
$M$, and we define the objects of $\soc_{\scrb}$ to be compact
$1$-manifolds (unions of circles and intervals) where each interval
has its boundary points labeled by an element of $\scrb.$ These labels
are meant to represent the boundary values of maps from these
intervals.
 
The $1$-morphisms in this category will be ``open-closed cobordisms.''
These are cobordisms between the compact manifolds with boundary that
respect the boundaries of the compact one-manifolds (in the sense that
the cobordisms restrict to give ordinary cobordisms of the boundaries
of the $1$-manifolds), as well as the labeling by elements of $\scrb$.
In other words, the boundary of an open-closed cobordism is
partitioned into three parts:
\begin{enumerate}
\item the ``incoming part,'' consisting of a number of circles and
  intervals, whose boundaries are labeled by elements of $\scrb$,
\item the``outgoing part'' of the same type, and,
\item the ``free boundary,'' which is itself a cobordism between the
  boundaries of the incoming and outgoing parts.  Each connected
  component of the free boundary is labeled by an element of $\scrb$,
  which respects the labeling of the boundaries of the incoming and
  outgoing parts.
\end{enumerate}
 
Finally, the $2$-morphisms are diffeomorphisms of the cobordisms that
respect all of this structure.
 
To state our main result, we intoduce some notation.  Let $F_{g,n}$ be
an oriented surface of genus $g$ with $n$ boundary components.  Let
${\bf m}$ be a fixed, finite set of points in $F_{g,n}.$ The set ${\bf
  m}$ is partitioned into subsets indexed by $\scrb$,
$$
{\bf m} = \bigcup_{b \in \scrb}m_b.
$$
We think of the subset $m_b$ as those points in the set ${\bf m}$
that are ``colored'' by the element $b \in \scrb$.  Let
$\diff(F_{g,n}, {\bf m}, \p)$ be the group of orientation-preserving
diffeomorphisms of $F_{g,n}$ that fix the boundary pointwise and
preserve each of the sets $m_b$ for $b \in \scrb.$ Let
$\Gamma_{g,n}^{\bf m}=\pi_0\diff(F_{g,n}, {\bf m}, \p)$ denote the
corresponding mapping class group.
 
Let $B\Gamma_{\infty, n}^{\scrb}$ be the colimit of the classifying
spaces, $\varinjlim_{g, {\bf m}} B\Gamma_{g,n}^{\bf m}.$ Here we are
taking the limit over the genus $g$ and the sets ${\bf m}$, which are
partially ordered by inclusion.  In this limit, the cardinality of
each subset $m_b$ tends to infinity.  The following is the main
theorem of this paper.
 
\begin{theorem}\label{main}  
  \begin{enumerate} 
  \item The category $\soc_{\scrb}$ of open and closed strings with
    set of $D$-branes $\scrb$ can be given the structure of a
    symmetric monoidal $2$-category, so its geometric realization
    $|\soc_{\scrb}|$ is an infinite loop space.
  \item For any $n$ there is a homology equivalence
    $$
    \bz \times \bz^{|\scrb|} \times B\Gamma_{\infty, n}^{\scrb}
    \xrightarrow{\cong_{H_*}}\Omega|\soc_{\scrb}|.
    $$ Here $|\scrb|$ is the cardinality of the set $\scrb$.  
  \item There is a homotopy equivalence of infinite loop spaces,
    $$
    \Omega|\soc_{\scrb}| \simeq  \Omega^\infty \bc \bp^\infty_{-1} \times \prod_{b \in \scrb} Q(\bc \bp^\infty_+)
    $$  where, as usual, $X_+$ denotes $X$
    with a disjoint basepoint and $Q(Y) =
    \varinjlim_{k}\Omega^k\Sigma^k(Y)$.
  \end{enumerate}
\end{theorem}

\med The paper is organized as follows.  In Section~1 we define the
notion of open-closed cobordism carefully and define the category
$\soc$ corresponding to the case when there is only one $D$-brane.
That is, the set $\scrb$ consists of one element.  This simplifying
assumption will allow us not to be concerned with the labeling of the
boundaries of the $1$-manifolds.  In this section we also show that
the category $\soc$ has a symmetric monoidal structure.  As the reader
will see, our constructions are direct generalizations of those of
Tillmann \cite{tillmann}, \cite{tillmann2}.  However the objects in
our category will have to be stratified in a rather delicate way, in
order to keep track of the boundary configurations in the cobordisms.
In Section~2 we adapt Tillmann's argument involving Harer stability to
prove part~2 of Theorem~\ref{main} when there is only one $D$-brane.
In this setting, part~3 of Theorem~\ref{main} will follow immediately
from part 2 using the Madsen-Weiss theorem \cite{madweiss}, and a
result of B{\"o}digheimer and Tillmann \cite{bodtill} which we review.
In Section~3 we describe the necessary modifications to prove
Theorem~1 for a general set $\scrb$ of $D$-branes.

% --------------------------------------------------------------------
\section{Open-closed cobordisms and the category $\soc.$}
\label{sec:catsoc}
% --------------------------------------------------------------------

We begin by defining the notion of an open-closed cobordism.

\begin{definition}
  A smooth (two-dimensional) \defem{open-closed cobordism}
  $\Sigma=(\Sigma,\bdry^\pm \Sigma, \bdryF \Sigma)$ is a smooth
  oriented surface $\Sigma,$ together with a choice of subsets
  $\bdry^-\Sigma,$ $\bdry^+\Sigma,$ $\bdryF \Sigma$ of its boundary
  $\bdry \Sigma,$ such that
  \begin{enumerate}
  \item $\bdry \Sigma = \bdry^- \Sigma\cup\bdry^+ \Sigma\cup\bdryF \Sigma,$ 
    
  \item each of $\bdry^-\Sigma, \bdry^+\Sigma, \bdryF \Sigma$ is an embedded
    one-dimensional manifold with boundary (to which we assign the
    orientation induced from $\bdry \Sigma$), and,
    
  \item $\bdryF \Sigma$ is a cobordism from $\bdry(\bdry^- \Sigma)$ to
    $\bdry(\bdry^+ \Sigma).$
  \end{enumerate}
  The sets $\bdry^-\Sigma,$ $\bdry^+\Sigma,$ $\bdryF\Sigma$ are called the
  incoming, outgoing, and free boundaries, respectively. We will write
  $\bdry^\pm\Sigma$ for $\bdry^-\Sigma\cup\bdry^+\Sigma.$
  
  Given two such cobordisms $\Sigma_1$ and $\Sigma_2,$ we say that
  they are \defem{isomorphic} if there exists an
  orientation-preserving diffeomorphism $\phi:\Sigma_1\to\Sigma_2$
  which restricts to diffeomorphisms
  $\bdry^-\Sigma_1\to\bdry^-\Sigma_2,$
  $\bdry^+\Sigma_1\to\bdry^+\Sigma_2,$
  $\bdryF\Sigma_1\to\bdryF\Sigma_2.$ Denote the space of such
  isomorphisms by $\diffocp(\Sigma_1,\Sigma_2).$ If we also have given
  diffeomorphisms $\bdry^-\Sigma_1 \cong \bdry^-\Sigma_2$ and
  $\bdry^+\Sigma_1 \cong \bdry^+\Sigma_2,$ we define
  $\diffocp(\Sigma_1,\Sigma_2;\bdry^\pm)$ to be the set of elements of
  $\diffocp(\Sigma_1,\Sigma_2)$ which restrict to the fixed
  diffeomorphisms on $\bdry^\pm\Sigma_1=\bdry^\pm\Sigma_2. $
\end{definition}

% --------------------------------------------------------------------
 
% --------------------------------------------------------------------

\med 
We now go about defining the category of open and closed strings.
Actually there will be several versions of this category which will
have different applications.

The objects in such a category are essentially diffeomorphism classes
compact $1$-manifolds, which may or may not have boundary.  Such a
manifold is diffeomorphic to a disjoint union of circles and
intervals.  However, as we will see below, the order of these
manifolds, and even various groupings of the intervals, will be
important.  So we will define a ``general object'' to be a sequence of
diffeomorphism classes of circles and intervals, which we enumerate by
sequences of zeros and ones (standing for the circles and intervals,
respectively).

\begin{definition}\label{general}
  A \defem{general object} is a sequence $\vec{n}=(n_1,\ldots,n_k),$
  where $n_i\in\set{0,1}.$ The length of the sequence $k$ can be an
  arbitrary nonnegative integer. Each such tuple $\vec{n}$ stands for
  a one-dimensional submanifold $C_{\vec{n}}=\bigcup_{i=1}^k I(i,n_i)$
  of $\R^2,$ where we define $I(i,0)$ to be the circle of radius
  $\nicefrac14$ centered at $(i,0)$ and
  $I(i,1)=[i-\nicefrac14,i+\nicefrac14]\cross\set{0}.$ We will use the
  notation $\ell(\vec{n})=k$ for the length of a tuple $\vec{n}$.
  
  A \defem{general morphism} between $\vec{n}$ and $\vec{m}$ is an
  open-closed cobordism $\Sigma$ smoothly embedded in $\R^3,$ with
  $\bdry^-\Sigma=C_{\vec{n}}\cross\set{0}$ and
  $\bdry^+\Sigma=C_{\vec{m}}\cross\set{t}$ for some $t>0.$ In
  addition, we require that:
  \begin{enumerate}
  \item $\bdry^\pm\Sigma$ has a collar, in such a way that for some
    $\mu>0$ we have $C_{\vec{n}}\cross[0,\mu),C_{\vec{m}}\cross(t-\mu,t]
    \subset\Sigma.$
    
  \item The orthogonal projection  in the $t$-direction is a Morse function
    on the interior of $\Sigma.$ 
    
  \item The same projection is a Morse function when restricted to
    $\bdryF\Sigma.$
  \end{enumerate}
\end{definition}

\begin{figure}[ht]
  \centering
  \includegraphics[height=6cm]{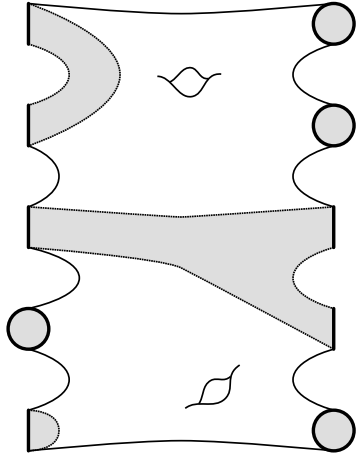}
  \caption{An open-closed cobordism from $\vn = (1,1,1,0,1)$ to $\vm=(0,0,1,1,0)$}
  \label{fig:figone}
\end{figure}

The first condition in the definition of general morphism allows us to
define compositions in the following way.  Given
$\Sigma_1\subset\R^2\cross[0,t_1]$ a general morphism from $ \vec{n}$
to $\vec{m}$,  and $\Sigma_2\subset\R^2\cross[0,t_2]$ a general morphism
from $ \vec{m}$ to $\vec{k},$ define their composition
$\Sigma_2\comp\Sigma_1\subseteq\R^2\cross[0,t_1+t_2]$ as the union of
$\Sigma_1$ with the result of translating $\Sigma_2$ by $t_1$ in
the $t$-direction. This clearly yields an associative composition.

The categories of open and closed strings we will consider will
actually be $2$-categories.  We therefore need the following
definition.

\begin{definition}\label{twomorph} Given $\Sigma_1$ and $\Sigma_2$ general morphisms from
  $\vec{n}$ to $\vec{m},$ we define the space of general
  $2$-morphisms, $\nmor{2}_{gen}(\Sigma_1,\Sigma_2)=\emptyset$ unless
  $\Sigma_1$ and $\Sigma_2$ are isomorphic, in which case we define
  $\nmor{2}_{gen}(\Sigma_1,\Sigma_2)$ to be the set of
  $\phi\in\diffocp(\Sigma_1,\Sigma_2;\bdry^\pm)$ which fix some collar
  of $\bdry^\pm\Sigma_1\isom\bdry^\pm\Sigma_2,$ where we identify
  collars of the target 1-manifolds of $\Sigma_1$ and $\Sigma_2$ by a
  suitable linear scaling in the $z$ direction.
\end{definition}

Notice that by composition of diffeomorphisms, for each pair of
general objects $\vec{n}$ and $\vec{m}$, we have a category
$\soc_{gen}(\vec{n}, \vec{m})$ of general morphisms and $2$-morphisms.
Furthermore, composition of morphisms defines functors
$$
\soc_{gen}(\vec{n}, \vec{m}) \times \soc_{gen}(\vec{m}, \vec{k}) \la \soc_{gen}(\vec{n}, \vec{k})
$$
which together define a $2$-category $\soc_{gen}$ with objects,
$1$-morphisms, and $2$-morphisms given by the general objects, general
morphisms, and general $2$-morphisms defined above.
 
\med We now define certain partitions associated to our general
objects.  We do this by considering a connected open-closed cobordism
(general morphism) $\Sigma$ from $\vec{n}$ to $\vo.$ Here the general
object $\vo$ is a sequence of length one and corresponds to a single
circle.
 
Let $\bdryI\Sigma$ be the set of those connected components of
$\bdry^\pm\Sigma$ which are intervals. In  the case we are considering  all of $\bdryI
\Sigma$ lies in the incoming part, $\bdry^-\Sigma$, and are in
bijection with the ``$1$'' entries in the sequence $\vec{n}.$ Notice
that the cobordism $\Sigma$ defines a partition of $\bdryI\Sigma$ into
subsets consisting of those intervals that lie in the same path
component of the boundary, $\p \Sigma$.

We describe this combinatorially as follows. 
\begin{definition}\label{permute}
  Let $k$ be the length of a general object $\vec{n}$. Let $I(\vec{n})
  = (i_1, \cdots, i_{\alpha (\vec{n})})$ be the ordered subset of
  $\{1, \cdots, k\}$ defined by $j \in I(\vec{n})$ if and only if $n_j
  = 1$.  So the length $\alpha (\vec{n})$ of $I(\vec{n})$ is the
  number of ones in the sequence $\vec{n}$.
  
  Consider the permutation $\sigma_\Sigma (\vec{n})$ of $I(\vec{n})$
  defined as follows.  Let $j \in I(\vec{n})$.  By definition this
  represents an interval $I_j \subset \bdryI\Sigma$.  $I_j$ lies in a
  connected component of the topological boundary $\bdry\Sigma.$
  Define $\sigma_\Sigma (j) = \ell$, where $I_{\ell}$ is the interval
  lying next to $I_j$ in this component, in the ordering induced by
  the orientation of $\Sigma.$ We call this permutation an \defem{open
    boundary permutation} of $\vec{n}$, or more precisely, of the
  one-manifold $C_{\vec{n}} = \p^+\Sigma$.  We view this permutation
  $\sigma_\Sigma$ as an element of the symmetric group $S_{\alpha
    (\vec{n})}$.
\end{definition}

\med

\begin{figure}[ht]
  \centering
  \includegraphics[height=6cm]{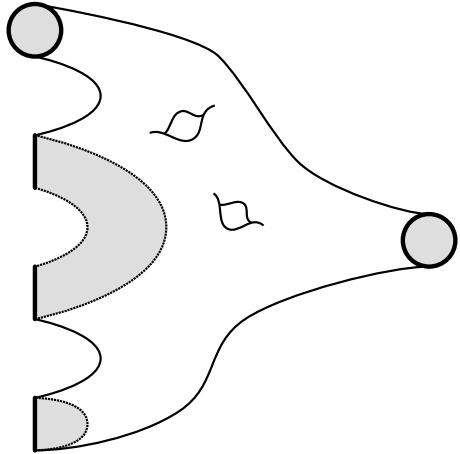}
  \caption{$\vn = (0,1,1,1), \, \sigma = (2,3)(4) \in S_{\alpha(\vn)}$}
  \label{fig:figtwo}
\end{figure}

\med Notice that different open-closed cobordisms $\Sigma$ from
$\vec{n}$ to $\vo$ can define the same permutation $\sigma \in
S_{\alpha (\vec{n})}$, and that every element of the symmetric group
can arise in this way.

\med We now define how open boundary permutations on $\vn$ may be
pulled back by a general morphism $\Sigma_{\vm,\vn}$ from $\vm$ to
$\vn.$
\begin{definition}
  Let $\Sigma_{\vm,\vn}$ be a general morphism from $\vm$ to $\vn$.
  Let $\sigma \in S_{\alpha (\vn)}$ be the open boundary permutation
  of $\vn$ and choose a connected general morphism $\Sigma_{\vn, \vo}$
  from $\vn$ to $\vo$ inducing $\sigma.$ Define the pullback
  permutation $\Sigma_{\vm, \vn}^*(\sigma) \in S_{\alpha (\vm)}$ to be
  the open boundary permutation of $\vm$ induced by the glued
  cobordism $\Sigma_{\vm, \vn}\circ \Sigma_{\vn, \vo}$ from $\vm$ to
  $\vo$.
\end{definition}
We leave it as an easy exercise for the reader to check that this
notion of pullback is well-defined.

\med Next, we define a 2-category $\socD$ as a refinement of
$\soc_{gen}$ that incorporates these open boundary permutations.

\begin{definition}\label{socD}
  Define the 2-category $\socD$ to have objects consisting of pairs
  $(\vn, \sigma)$, where $\vn$ is a general object and $\sigma \in
  S_{\alpha (\vn )}$.  Define the 1-morphisms, $\socD((\vec{n},
  \sigma), (\vec{m}, \tau))$ to be those general morphisms $\Sigma$
  from $\vn$ to $\vm$ such that $\Sigma^*(\sigma) = \tau$. The
  2-morphisms from $\Sigma_1$ to $\Sigma_2$ are the general
  2-morphisms, i.e., diffeomorphisms from $\Sigma_1$ to $\Sigma_2$ as
  in definition \ref{twomorph}.
\end{definition}

\med Define the category $\socG$ as the quotient $\socD/\!\!\equiv,$
where $\equiv$ is the equivalence relation that identifies pairs of
2-morphisms in the same path component in the corresponding
$\diffocp(\Sigma_1,\Sigma_2;\bdry^\pm).$ This category has the same
objects and morphisms as $\socD$, with each space
$\diffocp(\Sigma_1,\Sigma_2;\bdry^\pm)$ of 2-morphisms replaced by its
set of connected components, the mapping class groups which we denote
by $\Gamma(\Sigma_1,\Sigma_2;\bdry^\pm).$ We refer to \cite{tillmann2}
for a treatment of the relevant quotient construction on 2-categories.
The following result is immediate.

\med
\begin{proposition}
  The canonical map $\pi_0:\socD\to\socG$ is a 2-functor between
  2-categories.
\end{proposition}

% --------------------------------------------------------------------
%\subsection{The monoidal pairing}
%\label{sec:pairing}
% --------------------------------------------------------------------
Now we describe a pairing on $\socD$ which is not well-defined, nor
associative, nor functorial. However, this pairing will descend to a
symmetric monoidal structure on a subquotient of $\socG$.

\med
Given $\vec{n_1},\vec{n_2}\in\ob(\socD),$ define
$\vec{n_1}\tensor\vec{n_2}$ to be the juxtaposition of the two tuples
with $\vec{n_1}$ followed by $\vec{n_2}.$ Given
$\Sigma_1\in\socD(\vec{n_1},\vec{m_1}),$
$\Sigma_2\in\socD(\vec{n_2},\vec{m_2}),$ we construct
$\Sigma_1\tensor\Sigma_2\in\mor{\socD}(\vec{n_1}\tensor\vec{n_2},\vec{m_1}\tensor\vec{m_2})$
by taking the union of $\Sigma_1$ and $\Sigma_2$ after they have been
modified by suitable diffeomorphisms of $\R^3$ to make them be
disjoint and have the same height. More precisely, construct
$\Sigma_1\tensor\Sigma_2$ as follows. Suppose that $\Sigma_i\subset
\R^2\times[0,t_i].$ Choose a collar-preserving diffeomorphism
$r:[0,t_2]\to[0,t_1];$ in this way the diffeomorphism
$R(x,y,z):=(x,y,r(z))$ takes $\Sigma_2$ to an open-closed cobordism
$R(\Sigma_2)\in\socD(\vec{n_2},\vec{m_2})$ which lies in
$\R^2\times[0,t_1].$ Next, choose a smooth orientation-preserving
function $s:[0,t_1]\to\R_{\ge0}$ which equals $1+\ell(n_1)$ in a
neighborhood of $0$ and $1+\ell(m_1)$ in a neighborhood of $t_1,$ in
such a way that the map $S(x,y,z):=(x+s(z),y,z)$ takes $R(\Sigma_2)$
into an smooth open-closed cobordism $\tSigma_2$ with collar
which is disjoint from $\Sigma_1.$ Then, define
$\Sigma_1\tensor\Sigma_2:=\Sigma_1\cup \tSigma_2.$ As in the
closed case \cite{tillmann}  this yields a well-defined map of 2-morphisms which
descends to a map of 2-morphisms in $\socG.$

% --------------------------------------------------------------------
% \subsection{The quotient category}
% --------------------------------------------------------------------
To make this pairing into a strictly associative bifunctor, we will
make an identification among surfaces. As described for the closed
case \cite{tillmann2}, we will identify two surfaces in
$\socD(\vec{n},\vec{m})$ if they are related by a sequence of moves of
the same types as those described in Section~3 of \cite{tillmann2}. We
describe these moves next.

Let $\Sigma\subseteq\R^2\cross[0,t]$ be a smooth open-closed
cobordism, and pick $n$ planes $P_k=\R^2\cross\set{t_k}$ where
$0<t_1<\cdots<t_n<t,$ in such a way that each of the sets
$\Sigma_0,\ldots,\Sigma_n$ into which the planes divide $\Sigma$ is
itself a smooth open-closed cobordism between 1-manifolds $\Sigma\cap
P_{k}$ and $\Sigma\cap P_{k+1}.$ Here, $P_0 := \R^2\cross\set{0},
P_{n+1} := \R^2\cross\set{t}.$ Let the connected components of
$\Sigma_k$ be $\Sigma_{k,i},$ $i=1,\ldots,j_k.$ Allow each
$\Sigma_{k,i}$ to be independently transformed by a sequence of
operations of the following two types, as long as the resulting images
$\tSigma_{k,i}$ are pairwise disjoint:
\begin{itemize}
  \item \emph{$R$-move (Rescaling the height):} This is an operation of the
  form $R_{k,i}(x,y,z) = (x,y,r_{k,i}(z))$ on each connected component
  $\Sigma_{k,i},$ where
  $r_{k,i}:[t_k,t_{k+1}]\to[t_k,\tilde{t}_{k,i}]$ is a diffeomorphism
  with derivative one on a neighborhood of $\bdry[t_k,t_{k+1}].$ This
  rescales the height of all the connected components of $\Sigma_k$ to
  make it lie in $\R^2\cross[t_k,\tilde{t}_{k+1}];$ the slices
  $\Sigma_{k+1},\ldots,\Sigma_n$ are understood to be translated by
  $\tilde{t}_{k+1}-t_{k+1}$ in the $z$ direction.
  
\item \emph{$S$-move (Shifting parallel to the $xy$-plane):} This is
  an operation of the form $S_{k,i}(x,y,z) = (x,y,z) +
  (s^1_{k,i}(z),s^2_{k,i}(z),0),$ where the
  $s^j_{k,i}:[t_k,t_{k+1}]\to\R^2$ are smooth maps which vanish on a
  neighborhood of $\bdry[t_k,t_{k+1}].$ In particular, this move
  allows passing different connected components of $\Sigma$ across
  each other.
\end{itemize}
Define a 2-category $\socDb$ with the same objects as $\socD,$ by
letting $\socDb((\vec{m}, \sigma), (\vec{n}, \tau))$ be the full
subcategory of $\socD((\vec{m}, \sigma), (\vec{n}, \tau))$ consisting
of those open-closed cobordisms $\Sigma$ such that no connected
component of $\Sigma$ is an open-closed cobordism to the empty
1-manifold. Define the 2-category $\socGb$ analogously, and define
$[\socDb],$ $[\socGb]$ to be the quotient categories under the $R$ and
$S$ moves.

\begin{theorem} \label{symmonoidal}
  The monoidal pairing above  descends
  to $[\socGb]$ and makes it a strict symmetric  monoidal  
  2-category.
\end{theorem}

We omit the proof of this theorem, which is a formal argument using
the following lemma, proceeding exactly as the proof of Theorem~3.3 in
\cite{tillmann2}.
\begin{lemma}
  Let $\Sigma\in\socD((\vec{m}, \sigma), (\vec{n}, \tau)),$ and let
  $E(\Sigma)\subset \diffocp(\Sigma,\bdry^\pm)$ be the subgroup of
  automorphisms which are obtained as the composition of a sequence of
  $R$- and $S$-moves.  Then, $E(\Sigma)$ is contained in the connected
  component of the identity.
\end{lemma}
\begin{proof}
  The argument is the same as the one in $\cite{tillmann2}$ for the
  closed case.  One argues, from the form of the allowable moves, that
  for $\Sigma$ connected, any $\phi\in E(\Sigma)$ must fix each
  critical point of the $t$-axis projection, and hence that $\phi$
  fixes pointwise the subset $S$ of the union of the critical level
  sets which consists of those connected components containing a
  critical point. Since the complement of $S$ in $\Sigma$ is a union
  of cylinders over closed 1-manifolds with boundary, the result
  follows because the $R$- and $S$-moves cannot twist these cylinders.
  For $\Sigma$ not connected but lying in $\socDb,$ the argument goes
  through because the condition that no connected component of
  $\Sigma$ is a cobordism to the empty 1-manifold guarantees that
  $S$-moves cannot swap two connected components of $\Sigma.$
\end{proof}

\med As in \cite{tillmann2}, given a 2-category $\catC,$ we denote by
$\B\catC$ the category enriched over topological spaces in which
$\ob(\B\catC)=\ob(\catC)$ and $\B\catC(x,y)$ equals the realization of
the category $\catC(x,y).$ We will use $|\catC|$ to denote the
bisimplicial realization $B\B\catC$ of the nerve of $\B\catC.$

For the remainder of the paper, we ease the notation and let the
category of open and closed strings, $\soc$, be the the category
$\socGb$ above.  The following is an immediate consequence of
theorem~\ref{symmonoidal} using Segal's well known relationship
between symmetric monoidal structures on a category and infinite loop
structures on its classifying space \cite{segal3}.

\med
\begin{theorem}\label{infiniteloop} The classifying space of the category of open and closed strings,  $|\soc|$,  is an infinite loop space.
\end{theorem}

% --------------------------------------------------------------------
\section{The homology type of the category $\soc$}
\label{sec:hlgytype}
% --------------------------------------------------------------------
As in the introduction, let $F_{g,n}^m$ denote a fixed oriented
surface of genus $g$, with $n$ boundary components, and $m$ marked
points.  Let $\diff(F_{g,n}^{(m)}, \p)$ be the group of
orientation-preserving diffeomorphisms that fix the boundary
pointwise, and the marked points only setwise.  That is, the
diffeomorphisms may permute the marked points.  Let
$\Gamma_{g,n}^{(m)} = \pi_0(\diff(F_{g,n}^{(m)}, \p))$ be the
corresponding mapping class group.  We consider the colimits of the
classifying spaces as $g$ and $m$ increase,
 $$
B \Gamma_{\infty, n}^{(\infty)} = \varinjlim_{g, m \to \infty}B\Gamma_{g,n}^{(m)}.
$$
The main goal of this section is to prove the following two results.

\med
\begin{theorem}\label{homology}
For each $n$ there is a homology equivalence,
$$
\alpha : \bz \times \bz \times B \Gamma_{\infty, n}^{(\infty)} \xrightarrow{\cong_{H_*}} \Omega B\soc.
$$
\end{theorem}

\med
\begin{theorem}\label{homotopy}
For each $n$ there is a homotopy equivalence
$$
\bz \times \bz \times (B \Gamma_{\infty, n}^{(\infty)})^+ \simeq \Omega^\infty ( \bc \bp^\infty_{-1})  \times Q(\bc \bp^\infty_+)
$$
where the superscript ``$+$'' refers to the Quillen plus construction with respect to the maximal perfect subgroup of the fundamental group.  
\end{theorem}

\med As we will see below, theorem~\ref{homotopy} follows quickly from
a splitting theorem for stable mapping class groups due to
B\"odigheimer and Tillmann \cite{bodtill}, and Madsen and Weiss's
proof of the stable Mumford conjecture.  So most of this section will
be devoted to the proof of Theorem~\ref{homology}.

 \med
\begin{proof}[Proof of Theorem \ref{homology}.]   
  Consider the morphism category, $\socD((\vn, \sigma),\vo)$ Remember
  that the $1$-manifold $C_{\vo}$ is a single circle, so that there
  are no boundary permutations for $\vo.$ Notice that on the
  classifying space level we have
  $$
  B\socD((\vn, \sigma),\vo) \simeq \bigdisjcup_{\Sigma}
  B\diffocp(\Sigma;\bdry^\pm)
  $$
  where $\Sigma$ ranges over the isomorphism types of open-closed
  cobordisms from $C_{\vec{n}}$ to $C_{\vo}.$ that induce the boundary
  permutation $\sigma \in S_{\alpha (\vn)}$.  In the absence of open
  strings, the isomorphism type of a connected cobordism from $n$
  circles to $1$ circle is completely determined by its genus. In our
  present setting, there is an another invariant of connected
  open-closed cobordisms in addition to genus. This is the number $w$
  of ``windows'' of a an open-closed cobordism $\Sigma$, which is the
  number of connected components of the free boundary which are
  homeomorphic to circles.  See Figure~\ref{figthree}.

  \begin{figure}[ht]
    \centering
    \includegraphics[height=6cm]{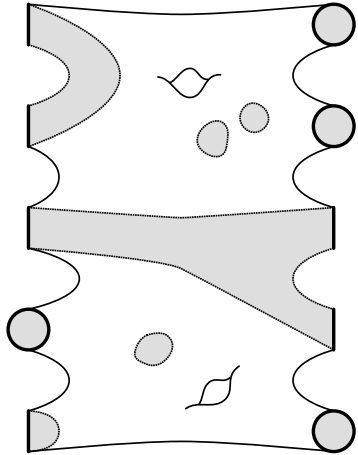}
    \caption{An open-closed cobordism with window number  $w = 3$.}
    \label{figthree}
  \end{figure}

We remark that if we had not fixed a boundary permutation $\sigma \in
S_{\alpha (\vn)}$, then the isomorphism type of a connected
open-closed cobordism $\Sigma$ from $\vn$ to $\vo$ would be determined
by \emph{three} invariants: its genus, the number of windows, and the
boundary permutation $\sigma (\Sigma) \in S_{\alpha (\vn)}$.
 
\med Let $\Sigma_{(\vm, \tau), (\vn, \sigma)}$ be an open-closed
cobordism representing a 1-morphism in $\socD$.  Let $\phi\in
\diffocp(\Sigma_{(\vm, \tau), (\vn, \sigma)}, \bdry^\pm).$ The
connected components of the topological boundary of $\Sigma_{(\vm,
  \tau), (\vn, \sigma)}$, and the action of $\phi$ on them, are given
as follows:
\begin{enumerate}
\item There are the circle components of the
  incoming and outgoing boundary. By definition, these are fixed
  pointwise by $\phi.$
\item There are $w$ windows, i.e., circle components of the free
  boundary. These are only fixed by $\phi$ as a set.  It follows that
  $\phi$ is isotopic to a diffeomorphism of the corresponding surface
  with punctures replacing the ``windows''.  By a ``puncture'' we mean
  a point that has been removed from the surface, where in the case
  of a ``window'', an open disc has been removed from the surface.
  Such a diffeomorphism may permute these punctures.
 
\item Call the components of the topological boundary $\p
  \Sigma_{(\vm, \tau), (\vn, \sigma)}$ that contain at least one
  incoming or outgoing interval ``mixed'' boundary components.  Since
  $\phi$ fixes each of the incoming and outgoing boundaries pointwise,
  it follows that $\phi$ restricts to a self-diffeomorphism on each of
  the mixed components.  Moreover, since each of these diffeomorphisms
  fixes a whole interval, $\phi$ is isotopic to a diffeomorphism that
  fixes each of the mixed components pointwise.
\end{enumerate}

From this description, we have that that 
$$B\pi_0 \diffocp(\Sigma_{(\vm, \tau), (\vn, \sigma)} ,\bdry^\pm) \htpyeq B\Gamma_{g,c(\Sigma)}^w,$$
where $c(\Sigma)$ is the number of boundary components of $\p \Sigma$
of the first or second type, and $w$ is the number of windows, or
punctures.  In particular for an open-closed cobordism from $(\vm,
\tau)$ to $\vo$, the number $c (\Sigma)$ is the number of zeros in the
sequence $\vm$ plus the number of disjoint cycles in $\tau \in
S_{\alpha (\vm)}$, plus one.  Since, in this case this number only
depends on the object $(\vm, \tau)$, we call it $c((\vm, \tau))$.
 
Thus for an object $(\vm, \tau)$ of the open-closed string category
$\soc$, we have proven the following.
 
\begin{lemma}\label{class}
  Let $\B\soc((\vm, \tau), \vo)$ be the classifying space of the
  category of morphisms from the object $(\vm, \tau)$ in $\soc$ to the
  object $\vo$.  Then
  $$
  \B\soc((\vm, \tau), \vo) \simeq \coprod_{w, \, g \geq 0}
  B\Gamma_{g, c(\vm, \tau)}^w.
  $$
\end{lemma}
 
\med
 
Fix an open-closed cobordism $T$ in $\soc(\vo, \vo)$ with one window
and genus one. See Figure~\ref{figfour}.

\begin{figure}[ht]
  \centering
     \includegraphics{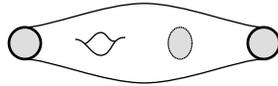}
  \caption{The surface $T$}
  \label{figfour}
\end{figure}

By gluing $T$ on the outgoing boundary and extending isomorphisms by
the identity, we obtain a functor from $\soc((\vec{n}, \sigma), \vo)$
to itself, and hence a map
\begin{equation}\label{tt}
t: \B\soc((\vec{n}, \sigma), \vo)\la \ \B\soc((\vec{n}, \sigma), \vo).
\end{equation}

Let $\B\soc_\infty((\vec{n}, \sigma))$ be the homotopy colimit
of the system
$$\B\soc((\vec{n}, \sigma), \vo)\stackrel{t}{\to} \B\soc((\vec{n},
\sigma), \vo)\stackrel{t}{\to}\cdots.$$
Since the maps in this colimit
increase both the genus and the number of windows, Lemma~\ref{tt}
implies the following.

\med
\begin{lemma}\label{colim}
  $$
  \B\soc_\infty((\vec{n}, \sigma)) \simeq \bz \times \bz  \times B\Gamma_{\infty, c(\vn, \sigma)}^{(\infty)}.
  $$
\end{lemma}
Notice that in the statement of this lemma, the two copies of $\bz$
reflect the genus and the number of windows.

\med This definition of extends to a contravariant functor
$\B\soc_\infty$ from $\B\soc$ to spaces, since, clearly, a morphism
$\psi$ from $(\vn, \sigma)$ to $(\vm, \tau)$ defines a pull back map
$$
\psi^* : \B\soc_\infty((\vec{m}, \tau)) \to \B\soc_\infty((\vec{n}, \sigma)),
$$
which in view of Lemma~\ref{colim} is homotopic to the gluing map
$$
\psi^* : \bz \times \bz \times B\Gamma_{\infty, c(\vm, \tau)}^{(\infty)} \to \bz \times \bz \times B\Gamma_{\infty, c(\vn, \sigma)}^{(\infty)}.
$$
But Harer's stability theorem, as adapted by B\"odigheimer and
Tillmann \cite{bodtill} to the setting of infinitely many punctures,
implies that each such gluing map is a homology isomorphism.  By the
generalized group completion theorem (\cite{mcduffseg},
\cite{tillmann}) this implies the following.

\med
\begin{lemma}\label{hlgyfib} Let $\hocolim\B\soc_\infty$ be the homotopy colimit of the functor
  $\B\soc_\infty$.  Then the natural projection map
$$
p : \hocolim \B\soc_\infty \to B\soc
$$
is a homology fibration, with fiber homology equivalent to $ \bz
\times \bz \times B\Gamma_{\infty, n}^{(\infty)} $ for any positive
integer $n$.
\end{lemma}

\med Notice that $\hocolim \B\soc_\infty $ is the homotopy colimit
under the geometric realization of the maps induced by the operation
$t$ of the homotopy colimit of the contravariant functor $\B\soc_0$
which assigns to an object $(\vm, \tau)$, the space of morphisms,
$\B\soc((\vm, \tau), \vo)$.  But the homotopy colimit of this
functor is contractible since it is the geometric realization of a
category with a terminal object (the Quillen over-category $\soc /
\vo$).  Thus $\hocolim \B\soc_\infty $ is the homotopy colimit of
contractible spaces, and hence is itself contractible.
(See~\cite{tillmann} for this type of argument.)

Thus the homotopy fiber of the map $p:\hocolim \B\soc_\infty \to
B\soc$ is homotopy equivalent to the loop space, $\Omega B\soc$.
Theorem~\ref{homology} then follows from Lemma~\ref{hlgyfib}.
\end{proof}

\med We now address Theorem~\ref{homotopy}.  The basis for this
theorem is the following result of B\"odigheimer and Tillmann
(\cite{bodtill}, Theorem~3.1), which makes a clever use of the Harer
stability theorem \cite{harer}.

\med
\begin{theorem}\label{bodtill} 
  There is a homotopy equivalence, 
  $$(B\Gamma_{\infty, n}^{(k)})^+ \simeq ( B\Gamma_\infty)^+ \times
  B(S_k\wr S^1),$$
  where $S_k$ is the symmetric group on $k$ letters
  and $S_k\wr S^1$ is the corresponding wreath product.
\end{theorem}

We remark that in this theorem $(B\Gamma_\infty)^+$ reflects
$(B\Gamma_{\infty, n})^+$ for any $n$, since as argued in
\cite{bodtill}, an easy consequence of Harer's stability theorem is
that the homotopy type of the plus construction of the classifying
space of the stable mapping class group does not depend on how many
(fixed) boundary components there are.

\med For the sake of completeness, we give a sketch of the proof of
this theorem.  We refer the reader to \cite{bodtill} for details.  The
idea is to consider the central extension
$$
\bz^n \to \Gamma_{\infty, n} \to \Gamma^n_{\infty}
$$
which gives rise to a map $\theta : \Gamma^n_\infty \to (GL^+(2,
\br))^n = B(\bz^n) \simeq (S^1)^n.$ Here $GL^+$ stands for matrices
having positive determinant.  The map $\theta$ can be viewed as taking
a representative diffeomorphism of a surface with $n$ marked points,
$\phi : F \xrightarrow{\cong}F$, and considering the derivatives,
$D\phi_{x_1} : T_{x_i} F \xrightarrow{\cong} T_{x_i}F$.  By fixing
framings at the marked points, these can be viewed as elements of
$GL^+(2, \br)$.

In the mapping class group $\Gamma^{(n)}_\infty$, the representative diffeomorphisms can permute the marked points among themselves, and therefore one has a map
$$
\theta_n :  \Gamma^{(n)}_\infty \to S_n\wr GL^+(2, \br)
$$
which yields a fibration sequence
$$
B\Gamma_{\infty, n} \to B\Gamma^{(n)}_\infty \xrightarrow{\theta_n}
B( S_n\wr GL^+(2, \br))\simeq B(S_n\wr S^1).
$$

By forgetting marked points, one has a forgetful map $\phi :
B\Gamma^{(n)}_\infty \to B\Gamma_\infty$ so that the composite,
$B\Gamma_{\infty, n} \to B\Gamma^{(n)}_\infty \xrightarrow{\phi}
B\Gamma_\infty$ is a homology isomorphism by Harer's stability theorem
\cite{harer}.  Therefore, $\phi$ induces a homology equivalence
between this fibration and the trivial fibration,
$$
\phi : B\Gamma^{(n)}_\infty \xrightarrow{\cong_{H_*}} B\Gamma_\infty \times B( S_n\wr S^1).
$$
Upon applying Quillen's plus construction, this homology
equivalence becomes a homotopy equivalence, and Theorem~\ref{bodtill}
follows.

\med Now the classifying space $B(S_n\wr S^1) \simeq ES_k \times_{S_k}
(BS^1)^k = ES_k \times_{S_k} (\bc \bp^\infty)^k$, and the classical
result of Barratt-Quillen-Priddy is that when one takes the colimit
over $k$, there is a homology equivalence
$$
\bz \times \varinjlim_{k \to \infty} ES_k \times_{S_k} (\bc \bp^\infty)^k \xrightarrow{\cong_{H_*}} Q (\bc \bp^\infty_+).
$$
Passing to the Quillen plus construction we therefore have the
following:

\med
\begin{corollary}\label{infinite}  There is a homotopy equivalence
  $$
  \bz \times (B\Gamma_{\infty, n}^{(\infty)})^+ \simeq (
  B\Gamma_\infty)^+ \times Q(\bc \bp^\infty_+).
$$
\end{corollary}

\med

Now recall Madsen and Weiss's theorem, proving a generalization of Mumford's conjecture on the cohomology of the stable mapping class group.

\med
\begin{theorem}\label{madweiss}  There is a homotopy equivalence 
$$
\alpha : \bz \times ( B\Gamma_\infty)^+   \xrightarrow{\simeq} \Omega^\infty(\bc \bp^\infty_{-1}).
$$
of infinite loop spaces.
\end{theorem}

\med Combining Theorem~\ref{homology}, Corollary~\ref{infinite}, and
Theorem~\ref{madweiss}, we obtain Theorem~\ref{homotopy}.

\section{Using more $D$-branes}

Observe that taken together, Theorems~\ref{homology}
and~\ref{homotopy} of the last section imply the main
Theorem~\ref{main} in the case that the the set of $D$-branes $\scrb$
consists of only one element.  In this section we consider the case
when we have a finite set of $D$-branes $\scrb$, and we describe the
necessary modifications of the arguments in the last section needed to
prove Theorem~\ref{main} in this general setting.

\med In the presence of a set $\scrb$ of $D$-branes, we define our
general objects to be disjoint unions of circles and intervals, but
now each interval has its boundary components labeled by $D$-branes.
That is, we have the following definition.

\med
\begin{definition}  
  Given a finite set $\scrb$, a $\scrb$-general object is a sequence
  $\vn = (n_1, \cdots , n_k)$ where each $n_i \in \{0,1\}$, together
  with a labeling function of the boundary of the intervals of
  $C_{\vn}$,
  $$\beta: \bdry C_{\vn} \to \scrb.$$
\end{definition}

\med We similarly define a $\scrb$-general morphism to be an open
closed cobordism that respects the labeling functions.

\med
\begin{definition}  
  A $\scrb$-general morphism from $(\vn, \beta_1)$ to $(\vm, \beta_2)$ is a
  general morphism $\Sigma$ from $\vn$ to $\vm$ (in the sense of
  Definition~\ref{general}), together with a labeling function of its
  free boundary components,
  $$
  \beta_\Sigma : \bdryF\Sigma \to \scrb
  $$
  (where $\beta_\Sigma$ is constant on each connected component)
  such that the free boundary $\bdryF\Sigma$ is a labeled cobordism
  from $\bdry C_{\vn} = \bdry(\bdry^-\Sigma)$ to $\bdry C_{\vm} =
  \bdry(\bdry^- \Sigma).$ By ``labeled cobordism'' we mean that the
  labeling functions are compatible.  That is,
  $$
  \beta_1 \sqcup \beta_2 :  \bdry C_{\vn} \sqcup \bdry C_{\vm} \to \scrb
  $$
  is equal to the restriction of $\beta_\Sigma$ to the boundary,
  $$
  \beta_{\Sigma} : \bdry (\bdryF \Sigma ) \to \scrb.
  $$
 \end{definition}
 
 \med As we did in section 1, we need to stratify the objects further
 by considering ``open boundary permutations.''  Namely, we observe
 that if $(\Sigma, \beta_\Sigma)$ is a $\scrb$-general morphism
 between $(\vn, \beta)$ and $\vo$, (i.e an open-closed,
 $\scrb$-labeled cobordism from $C_{\vn}$ to $\vo$) then, as in
 definition \ref{permute} above, $\Sigma$ defines an open boundary
 permutation $\sigma_\Sigma \in S_{\alpha (\vn)}$.
 
 This allows us to define the open-closed string 2-categories
 $(\socD)_{\scrb}$,  $(\socG)_\scrb$, and ultimately $\soc_\scrb$ just
 like in section one, but the objects will now be $\scrb$-general
 objects together with open-boundary permutations, $(\vn, \beta,
 \sigma)$.  The 1-morphisms from $(\vn, \beta, \sigma)$ to $(\vm,
 \beta_2, \tau)$ are $\scrb$-general morphisms $(\Sigma,
 \beta_{\Sigma})$ from $(\vn, \beta)$ to $(\vm, \beta_2)$ that pull
 back the open boundary permutations, $\Sigma^*(\sigma) = \tau$. The
 $2$-morphisms $(\Sigma,\beta_\Sigma)\to(\Sigma',\beta_{\Sigma'}$ are
 the isomorphisms $\phi:\Sigma\to\Sigma'$ of open-closed cobordisms as
 defined in Definition~\ref{twomorph}, with the additional requirement
 that they preserve the labeling functions:
 $$\beta_{\Sigma'} \comp \phi = \beta_\Sigma : \bdryF \Sigma \to
 \scrb. $$
 The proof of Theorem~\ref{symmonoidal} then gives the
 following generalization.
 
 \med
 \begin{theorem} $\soc_{\scrb}$ is a strict symmetric monoidal $2$-category.
 \end{theorem}
 
 \med In the presence of the set $\scrb$ of $D$-branes,
 Theorem~\ref{homology} needs the following modifications.  Notice
 first of all that the windows in a $\scrb$-labeled open-closed
 cobordism $(\Sigma, \beta_\Sigma)$ (a $1$-morphism in $\soc_\scrb$)
 are connected components of the free boundary, and therefore are
 equipped with labels in $\scrb$:
 \begin{equation}\label{windows}
 \beta_{\Sigma} : W \to \scrb
 \end{equation}
 where $W \subset \pi_0(\bdryF(\Sigma))$ is the set of free boundary
 components which are windows.  The set $W$ is partitioned as
 $$
 W = \sqcup_{b \in \scrb} W_b
 $$
 where $W_b = \beta_\Sigma^{-1}(b)$.  Let $w_b \in \bz^+$ be the
 cardinality of the set $W_b.$ The analogue of Lemma~\ref{class} above
 is the following:
 
 \begin{lemma}\label{classB} 
   Let $\B\soc_\scrb((\vm, \beta, \tau), \vo)$ be the classifying
   space of the category of morphisms from the object $(\vm, \beta,
   \tau)$ in $\soc$ to the object $\vo$.  Then
   $$
   \B\soc((\vm, \beta,  \tau), \vo)  \simeq \coprod_{w = \{w_b\} \in (\bz^+)^{|\scrb |} }\quad \coprod_{ g \geq 0} B\Gamma_{g, c(\vm, \beta,\tau)}^w.
   $$
   Here the number $c(\vm, \beta, \tau)$ is defined completely
   analogously to the description of $c(\vm, \tau)$ prior to the
   statement of Lemma~\ref{class}.  The group $\Gamma_{g, c(\vm, \beta,
     \tau)}^w$ is the mapping class group of isotopy classes of
   diffeomorphisms of a surface of genus $g$, with $c(\vm, \beta, \tau)$
   boundary components, and $w = \Sigma_{b \in \scrb} w_b$ marked
   points, where $w_b$ of these marked points are labeled by the
   element $b$.  These diffeomorphisms can permute these marked
   points, but must preserve the labelings.
\end{lemma}

\med

\begin{figure}[ht]
  \centering
     \includegraphics{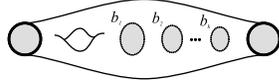}
  \caption{$T_\scrb$. The genus is one, and there is one window labeled by each $b \in \scrb$} 
  \label{figfive}
\end{figure}

One then defines a contravariant functor $(\B\soc_\scrb)_\infty$ from
$\B\soc_\scrb$ to spaces like what was done in Section 2.  However in
the presence of the set of $D$-branes $\scrb$ we define
$\B(\soc_\scrb)_\infty$ to be the homotopy colimit of the map
$$
t_\scrb : \B\soc((\vm, \beta, \tau), \vo) \to \B\soc((\vm, \beta,\tau), \vo)
$$
defined by gluing on the surface $(T_\scrb, \beta_{T_\scrb})$
defined to be a morphism in $\soc_\scrb (\vo, \vo) $ of genus one,
with exactly one window $W_b$ labeled by each element $b \in \scrb$.
See Figure~\ref{figfive}.

The remaining part of the proof of Theorem~\ref{homology} above goes
through to prove the following analogue:

\begin{theorem}\label{Bhomology} For each $n$ there is a homology equivalence
  $$
  \bz \times \bz^{|\scrb|} \times B\Gamma_{\infty, n}^{\scrb}
  \xrightarrow{\cong_{H_*}}\Omega|\soc_{\scrb}|.
 $$ where $ B\Gamma_{\infty, n}^{\scrb}$ is defined as in the statement of Theorem~\ref{main}.
\end{theorem}
 
\med Now just like Theorem~\ref{homotopy} above is a consequence of
B\"odigheimer and Tillmann's splitting of the classifying spaces of
stable mapping class groups (Theorem~3.1 of \cite{bodtill}) and Madsen
and Weiss's theorem (Theorem~\ref{madweiss} above) the following
analogue also follows from these theorems.
 
\med
\begin{theorem}\label{Bhomotopy}
  For each $n$ there is a homotopy equivalence
  $$
  \bz \times \bz^{|\scrb|} \times (B\Gamma_{\infty, n}^{\scrb})^+
  \simeq \Omega^\infty ( \bc \bp^\infty_{-1}) \times \prod_{b \in
    \scrb} Q(\bc \bp^\infty_+)
  $$
 \end{theorem}
 
 \med Finally we notice that Theorems~\ref{Bhomology}
 and~\ref{Bhomotopy} together are equivalent to the statement of
 Theorem~\ref{main}.

\end{document}